\documentclass[12pt,a4paper]{article}

\usepackage{amsmath}
\usepackage{float}
\usepackage{amsfonts}
\usepackage{amstext}
\usepackage{graphics}
\usepackage{epsf}
\usepackage[dvips]{epsfig}

\newtheorem{theorem}{Theorem}
\newtheorem{lemma}[theorem]{Lemma}
\newtheorem{corollary}[theorem]{Corollary}

\def \mod#1 #2{#1\ ({\rm mod}\ #2)}
\def \nodiv{{|\kern-3.9pt/}}

\def \prend{\vrule depth-1pt height7pt width6pt}

\def \endpf{{\ \ \prend \medbreak}}

\def \En{{\mathbb N}}

\title{Cubefree Binary Words Avoiding Long Squares}

\author{Narad Rampersad, Jeffrey Shallit, and Ming-wei Wang \\
Department of Computer Science \\
University of Waterloo \\
Waterloo, ON, N2L 3G1 \\
CANADA \\
{\tt nrampersad@math.uwaterloo.ca}\\
{\tt shallit@graceland.uwaterloo.ca}\\
{\tt m2wang@math.uwaterloo.ca} 
}

\begin{document}
\date{}
\maketitle

\centerline{\bf PLEASE NOTE:  AFTER THIS PAPER WAS PREPARED, WE LEARNED}

\centerline{\bf THAT ALL OUR RESULTS APPEARED (ALBEIT WITH DIFFERENT PROOFS)}

\centerline{\bf IN A PAPER OF F. M. Dekking, On repetitions of blocks in binary}

\centerline{\bf sequences, {\it J. Combin. Theory Ser. A} {\bf 20} (1976), 292--299.}

\begin{abstract}
Entringer, Jackson, and Schatz conjectured in 1974 that every infinite
cubefree binary word contains arbitrarily long squares.  In this paper
we show this conjecture is false:  there exist infinite cubefree binary
words avoiding all squares $xx$ with $|x| \geq 4$, and the number $4$
is best possible.  However, the Entringer-Jackson-Schatz conjecture is true
if ``cubefree'' is replaced with ``overlap-free''.
\end{abstract}

\section{Introduction}

Let $\Sigma$ be a finite nonempty set, called an {\it alphabet}.  We
consider finite and infinite words over $\Sigma$.  The set of all
finite words is denoted by $\Sigma^*$. The set of all infinite
words (that is, maps from $\En$ to $\Sigma$) is denoted by
$\Sigma^\omega$.  

     A {\it morphism} is a map $h: \Sigma^* \rightarrow \Delta^*$
such that $h(xy) = h(x)h(y)$ for all $x, y \in \Sigma^*$.  A morphism
may be specified by providing the image words $h(a)$ for all $a \in \Sigma$.
If $h:\Sigma^* \rightarrow \Sigma^*$ and
$h(a) = ax$ for some letter $a \in \Sigma$, then we say that
$h$ is {\it prolongable} on $a$, and we can then iterate $h$ infinitely
often to get the fixed point
$h^\omega(a) := a \, x \, h(x) \, h^2(x) \, h^3(x) \cdots $.

     A {\it square} is a nonempty word of the form $xx$, as in the
English word {\tt murmur}.  A {\it cube} is a nonempty word of the
form $xxx$, as in the English sort-of-word {\tt shshsh}.  An
{\it overlap} is a word of the form $axaxa$, where $x$ is a possibly
empty word and $a$ is a single letter, as in the English word
{\tt alfalfa}.

      It is well-known and easily proved that every word of length
$4$ or more over a two-letter alphabet contains a square as a subword.
However, Thue proved in 1906 \cite{Thue:1906} that there exist
infinite words over a three-letter alphabet that contain no squares;
such words are said to {\it avoid} squares or be {\it squarefree}.
Thue also proved that the word $\mu^\omega (0) = 0110100110010110\cdots$
is overlap-free (and hence cubefree); here
$\mu$ is the morphism sending $0 \rightarrow 01$ and $1 \rightarrow 10$.

      Entringer, Jackson, and Schatz \cite{Entringer&Jackson&Schatz:1974}
proved that while squares cannot be avoided over a two-letter
alphabet, arbitrarily long squares can.  More precisely, they proved
that there exist infinite binary words with no squares of length
$\geq 3$, and that the number $3$ is best possible.  Later, this
result was improved by Fraenkel and Simpson \cite{Fraenkel&Simpson:1995},
who proved that there exist
infinite binary words where the only squares are $00$, $11$, and $0101$.

      Entringer, Jackson, and Schatz conjectured in 1974 that
any infinite cubefree word over $\lbrace 0,1 \rbrace$
contains arbitrarily long
squares \cite[Conjecture B, p.\ 163]{Entringer&Jackson&Schatz:1974}.
In this paper we show that this
conjecture is false; there exist infinite cubefree binary words
with no squares $xx$ with $|x| \geq 4$.  The number $4$ is best
possible.  Further, we show that the Entringer-Jackson-Schatz
conjecture is true if the word ``cubefree'' is replaced with
``overlap-free''.

\section{A cubefree word without arbitrarily long squares}

     In this section we disprove the conjecture of Entringer, Jackson, and
Schatz.  First we prove the following result. 

\begin{theorem}
There is a squarefree infinite word over $\lbrace 0, 1, 2, 3 \rbrace$ with
no occurrences of the subwords $12$, $13$, $21$, $32$, $231$, or $10302$.
\label{first}
\end{theorem}

\begin{proof}
Let the morphism $h$ be defined by
\begin{eqnarray*}
0 &\rightarrow& 0310201023 \\
1 &\rightarrow& 0310230102 \\
2 &\rightarrow& 0201031023 \\
3 &\rightarrow& 0203010201
\end{eqnarray*}
Then we claim the fixed point $h^\omega(0)$ has the desired properties.

First, we claim that if $w \in
\lbrace 0, 1, 2, 3 \rbrace^*$ then $h(w)$ has no occurrences of 
$12$, $13$, $21$, $32$, $231$, or $10302$.  For if any of these
words occur as subwords of $h(w)$, they must occur within some $h(a)$ or
straddling the
boundary between $h(a)$ and $h(b)$, for some single letters $a, b$.
They do not;  this easy verification is left to the reader.

Next, we prove that if $w$ is any
squarefree word
over $\lbrace 0, 1, 2, 3 \rbrace$ having no occurrences of 
$12$, $13$, $21$, or $32$, then $h(w)$ is squarefree.

We argue by contradiction.  Let $w = a_1 a_2 \cdots a_n$ be a
squarefree string such that $h(w)$ contains a square,
i.e., $h(w) = xyyz$ for some $x, z \in \lbrace 0, 1, 2, 3 \rbrace^*$,
$y \in \lbrace 0, 1, 2, 3 \rbrace^+$.
Without loss of generality, assume that $w$ is a shortest such
string, so that $0 \leq |x|, |z| < 10$.  

Case 1:  $|y| \leq 20$.    In this case we can take $|w| \leq 5$.
To verify that $h(w)$ is squarefree,
it therefore suffices to check each of the 49 possible words $w \in
\lbrace 0, 1, 2, 3 \rbrace^5$ to ensure that $h(w)$ is squarefree in each case.

Case 2: $|y| > 20$.  First, we establish the following result.

\begin{lemma}
\begin{itemize}
\item[(a)]
Suppose $h(ab) = t h(c) u$ for some
letters $a, b, c \in \lbrace 0, 1, 2, 3 \rbrace$
and strings $t, u \in \lbrace 0, 1, 2, 3 \rbrace^*$.
Then this inclusion is trivial (that is,
$t = \epsilon$ or $u = \epsilon$) or $u$ is not a prefix
of $h(d)$ for any $d \in \lbrace 0, 1, 2, 3 \rbrace$.

\item[(b)]
Suppose there exist letters $a, b, c$ and
strings $s, t, u, v$ such that $h(a) = st$, $h(b) = uv$,
and $h(c) = sv$.  Then either $a = c$ or $b = c$.
\end{itemize}
\label{ming}
\end{lemma}

\begin{proof}
\begin{itemize}
\item[(a)]
This can be verified with a short computation.  In fact, the 
only $a, b, c$ for which the equality $h(ab) = t h(c) u$ 
holds nontrivially is $h(31) = t h(2) u$, and in this case
$t = 020301$, $u = 0102$, so $u$ is not a prefix of any $h(d)$.

\item[(b)]  This can also be verified with a short computation.
If $|s| \geq 6$, then no two distinct
letters share a prefix of length $6$.
If $|s|\leq 5$, then $|t| \geq 5$, and no two distinct letters 
share a suffix of length $5$.
\end{itemize}
\endpf
\end{proof}

     For $i = 1, 2, \ldots, n$ define $A_i = h(a_i)$.  
Then if $h(w) = xyyz$, we can write
$$h(w) = A_1 A_2 \cdots A_n = A'_1 A''_1 A_2 \cdots A_{j-1} 
A'_j A''_j A_{j+1} \cdots A_{n-1} A'_n A''_n$$ 
where
\begin{eqnarray*}
A_1 &=& A'_1 A''_1 \\
A_j &=& A'_j A''_j \\
A_n &=& A'_n A''_n \\
x &=& A'_1 \\
y &=& A''_1 A_2 \cdots A_{j-1} A'_j = A''_j A_{j+1} \cdots A_{n-1} A'_n \\
z &=& A''_n, \\
\end{eqnarray*}
where $|A''_1|, |A''_j| > 0$.
See Figure~\ref{fig1}.

\begin{figure}[H]
\begin{center}
\input cube1.pstex_t
\end{center}
\caption{The string $xyyz$ within $h(w)$ \protect\label{fig1}}
\end{figure}

     If $|A''_1| > |A''_j|$, then $A_{j+1} = h(a_{j+1})$ is a subword
of $A''_1 A_2$, hence a subword of $A_1 A_2 = h(a_1 a_2)$.  Thus
we can write $A_{j+2} = A'_{j+2} A''_{j+2}$ with
$$ A''_1 A_2 = A''_j A_{j+1} A'_{j+2}.$$ 
See Figure~\ref{fig2}.

\begin{figure}[H]
\begin{center}
\input cube2.pstex_t
\end{center}
\caption{The case $|A''_1| > |A''_j|$ \protect\label{fig2}}
\end{figure}

But then,
by Lemma~\ref{ming} (a), either $|A''_j| = 0$,
or $|A''_1| = |A''_j|$, or $A'_{j+2}$ is a not
a prefix of any $h(d)$.  All three conclusions are impossible.

     If $|A''_1| < |A''_j|$, then $A_2 = h(a_2)$ is a subword of
$A''_j A_{j+1}$, hence a subword of $A_j A_{j+1} = h(a_j a_{j+1})$.
Thus we can write $A_3 = A'_3 A''_3$ with
$$ A''_1 A_2 A'_3 = A''_j A_{j+1} .$$  
See Figure~\ref{fig3}.

\begin{figure}[H]
\begin{center}
\input cube3.pstex_t
\end{center}
\caption{The case $|A''_1| < |A''_j|$ \protect\label{fig3}}
\end{figure}

By Lemma~\ref{ming} (a), either $|A''_1| = 0$ or $|A''_1| = |A''_j|$
or $A'_3$ is not a prefix of any $h(d)$.  Again, all three conclusions
are impossible.

     Therefore $|A''_1| = |A''_j|$.  
Hence $A''_1 = A''_j$, $A_2 = A_{j+1}$, $\ldots$, $A_{j-1} = A_{n-1}$,
and $A'_j = A'_n$.  Since $h$ is injective, we have
$a_2 = a_{j+1}, \ldots, a_{j-1} = a_{n-1}$.
It also follows that $|y|$ is divisible by $10$ and
$A_j = A'_j A''_j = A'_n A''_1$.   But by Lemma~\ref{ming} (b), either
(1) $a_j = a_n$ or (2) $a_j = a_1$.  In the first case,
$a_2 \cdots a_{j-1} a_j = a_{j+1} \cdots a_{n-1} a_n$, so
$w$ contains the square $(a_2 \cdots a_{j-1} a_j)^2$, a contradiction.  In the
second case, $a_1 \cdots a_{j-1} = a_j a_{j+1} \cdots a_{n-1}$, so
$w$ contains the square $(a_1 \cdots a_{j-1})^2$, a contradiction.

     It now follows that
the infinite word 
$$h^\omega(0) = 0310201023 0203010201 0310230102 0310201023 0201031023  \cdots$$
is squarefree and contains no occurrences of $12$, $13$, $21$, $32$,
$231$, or $10302$.
\endpf
\end{proof}

\begin{theorem}
Let $\bf w$ be any infinite word satisfying the conditions of Theorem~\ref{first}.
Define a morphism $g$ by
\begin{eqnarray*}
0 &\rightarrow& 010011 \\
1 &\rightarrow& 010110 \\
2 &\rightarrow& 011001 \\
3 &\rightarrow& 011010
\end{eqnarray*}
Then $g({\bf w})$ is a cubefree word containing no squares $xx$ with
$|x| \geq 4$.
\label{second}
\end{theorem}

     Before we begin the proof, we remark that all the words
$12$, $13$, $21$, $32$, $231$, $10302$ must indeed be avoided,
because
\begin{eqnarray*}
g(12) &&\mbox{contains the squares}\ (0110)^2,\ (1100)^2,\ (1001)^2 \\
g(13) &&\mbox{contains the square}\ (0110)^2 \\
g(21) &&\mbox{contains the cube}\ (01)^3 \\
g(32) &&\mbox{contains the square}\ (1001)^2 \\
g(231)  &&\mbox{contains the square}\ (10010110)^2 \\
g(10302) &&\mbox{contains the square}\ (100100110110)^2.
\end{eqnarray*}

\begin{proof}
     The proof parallels the proof of Theorem~\ref{first}.   
     Let $w = a_1
a_2 \cdots a_n$ be a squarefree string, with no occurrences of 
$12$, $13$, $21$, $32$, $231$, or $10302$.  We first establish that if
$g(w) = xyyz$ for some
$x, z \in \lbrace 0, 1,2, 3 \rbrace^*$, $y \in \lbrace 0, 1, 2, 3 \rbrace^+$,
then $|y| \leq 3$.
Without loss of generality, assume $w$ is a shortest such string, so
$0 \leq |x|, |z| < 6$.

     Case 1:  $|y| \leq 12$.  In this case we can take $|w| \leq 5$.  To verify
that $g(w)$ contains no squares $yy$ with
$|y| \geq 4$, it suffices
to check each of the $41$ possible words $w \in \lbrace 0, 1,2, 3 \rbrace^5$.

      Case 2:  $|y| > 12$.  First, we establish the analogue of
Lemma~\ref{ming}.

\begin{lemma}
\begin{itemize}
\item[(a)]
Suppose $g(ab) = t g(c) u$ for some
letters $a, b, c \in \lbrace 0, 1, 2, 3 \rbrace$
and strings $t, u \in \lbrace 0, 1, 2, 3 \rbrace^*$.
Then this inclusion is trivial (that is,
$t = \epsilon$ or $u = \epsilon$) or $u$ is not a prefix
of $g(d)$ for any $d \in \lbrace 0, 1, 2, 3 \rbrace$.

\item[(b)]
Suppose there exist letters $a, b, c$ and
strings $s, t, u, v$ such that $g(a) = st$, $g(b) = uv$,
and $g(c) = sv$.  Then either $a = c$ or $b = c$, or
$a = 2$, $b = 1$, $c = 3$, $s = 0110$, $t = 01$, $u = 0101$, $v = 10$.
\end{itemize}
\label{ming2}
\end{lemma}

\begin{proof}
\begin{itemize}
\item[(a)]  This can be verified with a short computation.  
The only $a, b, c$ for which $g(ab) = t g(c) u$ holds nontrivially
are 
\begin{eqnarray*}
g(01) &=& 010 \ g(3) \ 110 \\
g(10) &=& 01 \ g(2) \ 0011 \\
g(23) &=& 0110 \ g(1) \ 10 .\\
\end{eqnarray*}
But none of $110$, $0011$, $10$ are prefixes of any $g(d)$.

\item[(b)] If $|s| \geq 5$ then no two distinct letters share a prefix
of length $5$.  If $|s| \leq 3$ then $|t| \geq 3$, and no two distinct
letters share a suffix of length $3$.  Hence $|s| = 4$, $|t| = 2$.
But only $g(2)$ and $g(3)$ share a prefix of length $4$, and only
$g(1)$ and $g(3)$ share a suffix of length $2$.
\end{itemize}
\endpf
\end{proof}

     The rest of the proof is exactly parallel to the proof
of Theorem~\ref{first}, with the following exception.  When we
get to the final case, where $|y|$ is divisible by $6$, we can use
Lemma~\ref{ming2} to
rule out every case except where $x = 0101$, $z = 01$,
$a_1 = 1$, $a_j = 3$, and $a_n = 2$.  Thus $w = 1 \alpha 3 \alpha 2$
for some string $\alpha \in \lbrace 0,1,2,3 \rbrace^*$.  This
special case is ruled out by the following lemma:

\begin{lemma}
     Suppose $\alpha \in \lbrace 0,1,2,3 \rbrace^*$, and let
$w = 1 \alpha 3 \alpha 2$.  Then either $w$ contains a square, or
$w$ contains an occurrence of one of the subwords $12$, $13$,
$21$, $32$, $231$, or $10302$.
\end{lemma}

\begin{proof}
      This can be verified by checking (a) all strings $w$ 
with $|w| \leq 4$, and (b) all strings of the form $w = abc w' de$, where
$a, b, c, d, e \in \lbrace 0, 1, 2, 3 \rbrace$ and $w' \in
\lbrace 0, 1, 2, 3 \rbrace^*$.  (Here $w'$ may be treated as an
indeterminate.)
\endpf
\end{proof}

    It now remains to show that if $w$ is squarefree and contains
no occurrence of $12$, $13$, $21$, $32$, $231$, or $10302$, then
$g(w)$ is cubefree.  If $g(w)$ contains 
a cube $yyy$, then it contains a square $yy$, and from what precedes
we know $|y| \leq 3$.  It therefore suffices to show that $g(w)$
contains no occurence of $0^3$, $1^3$, $(01)^3$, $(10)^3$,
$(001)^3$, $(010)^3$, $(011)^3$, $(100)^3$, $(101)^3$, $(110)^3$.
The longest such string is of length $9$, so it suffices to examine the
$16$ possibilities for $g(w)$ where $|w| = 3$.  This is left to the reader.

    The proof of Theorem~\ref{second} is now complete.
\endpf
\end{proof}

\begin{corollary}
    If $g$ and $h$ are defined as above, then
$$g(h^\omega(0)) = 
010011011010010110010011011001010011010110010011011001011010
\cdots$$
is cubefree, and avoids all squares $xx$ with $|x| \geq 4$.
\label{jcor}
\end{corollary}

\section{The constant $4$ is best possible}

     It is natural to wonder if the constant $4$ in Corollary~\ref{jcor}
can be improved.  It cannot, as the following theorem shows.

\begin{theorem}
Every binary word of length $\geq 30$ contains a cube
or a square $xx$ with $|x| \geq 3$.
\end{theorem}

\begin{proof}
This may be proved purely mechanically.  More generally,
let $P \subset \Sigma^*$ be a set of subwords to be avoided.
We create and traverse
a certain tree $T$, as follows.  The root of the tree is labeled
$\epsilon$.  If a node is labeled $x$ and contains no subword in $P$,
then it has children labeled $xa$ for each $a \in \Sigma$; otherwise it is
a leaf of $T$.    This tree is infinite if and only if there is an
infinite word avoiding the elements of $P$.

If $T$ is finite, then the height of $T$ gives the
length $l$ such that every word of length $l$ or greater contains
an element of $P$.  The tree can be created and
traversed using a queue and breadth-first search.

If the set $P$ is symmetric under renaming of the letters---as it is in this
case---we may further improve the procedure by labeling the root with any particular
letter, say $0$.  When we run this procedure on the statement of the
theorem, we obtain a tree with 289 leaves, the longest being
of length $30$.   The unique string of length $29$ starting with $0$ and
avoiding cubes and squares $xx$ with $|x| \geq 3$
is $00110010100110101100101001100$.
\endpf
\end{proof}

\section{Overlap-free words contain arbitrarily long squares}

     It is also natural to wonder if a result like Corollary~\ref{jcor}
holds if ``cubefree'' is replaced with ``overlap-free''.  It does not,
as the following result shows.

\begin{theorem}
Any infinite overlap-free word over $\lbrace 0,1 \rbrace$ contains
arbitrarily long squares.
\end{theorem}

\begin{proof}
By \cite[Lemma 3]{Allouche&Currie&Shallit:1998} we know that
if $\bf x$
is an overlap-free infinite word over $\lbrace 0, 1 \rbrace$, then there exist
a word $u \in \lbrace \epsilon, 0, 1, 00, 11 \rbrace$
and an overlap-free infinite word $\bf y$ such that
${\bf x} = u \mu({\bf y})$, where $\mu$ is the Thue-Morse morphism.  By iterating
this theorem, we get that every overlap-free infinite word must
contain $\mu^n(0)$ for arbitrarily large $n$; hence contains
arbitrarily long squares.
\endpf
\end{proof}

\section{Acknowledgments}

    We thank Jean-Paul Allouche for helpful discussions.

\newcommand{\noopsort}[1]{} \newcommand{\singleletter}[1]{#1}

\end{document}